\documentclass[12pt]{amsart}
\usepackage{amssymb,enumerate}
\usepackage{tikz}

\usepackage{pst-plot}
\usepackage{pst-node,pst-text,pst-3d,pstricks, hyperref}

\def\proof{\noindent{\bf{Proof.} }}
\def\sqr#1#2{{\vcenter{\hrule height.#2pt
        \hbox{\vrule width.#2pt height#1pt \kern#1pt
                \vrule width.#2pt}
        \hrule height.#2pt}}}
\def\square{\mathchoice\sqr64\sqr64\sqr{4}3\sqr{3}3}
\def\QED{\hfill$\square$\\}

\newtheorem{theorem}{Theorem}[section]

\newtheorem{corollary}[theorem]{Corollary}
\newtheorem{lemma}[theorem]{Lemma}
\newtheorem{proposition}[theorem]{Proposition}

\newtheorem{definition}[theorem]{Definition}
\newtheorem{remark}[theorem]{Remark}

\newtheorem{example}[theorem]{Example}
\newtheorem{question}[theorem]{Question}

\newtheorem{notation}[theorem]{Notation}

\newcommand{\m}{\mathfrak{m}}

\newcommand{\depth}{\ensuremath{{\rm{depth}}\;}}

\numberwithin{equation}{section}

\newcommand{\assg}[1]{\ensuremath{{\rm gr}_{R}{(#1)}}}

\def\dim{\mbox{\rm {dim}}}

\def\Ass{\mbox{\rm {Ass\,}}}

\def\max{\mbox{\rm {max\,}}}
\newcommand{\rar}{\rightarrow}

\newcommand{\sdepth}{\ensuremath{{\rm{sdepth}\;}}}
\begin{document}
\title{A Lower Bound For Depths of Powers of Edge Ideals}

\author{Louiza Fouli}
\address{Department of Mathematical Sciences \\
New Mexico State University\\
P.O. Box 30001 \\
Las Cruces, NM 88003}
\email{lfouli@math.nmsu.edu}
\urladdr{http://www.math.nmsu.edu/~lfouli}

\author{Susan Morey}
\address{Department of Mathematics \\
Texas State University\\
601 University Drive\\
San Marcos, TX 78666}
\email{morey@txstate.edu}
\urladdr{http://www.txstate.edu/~sm26/}


\thanks{The first author was partially supported by a grant from the Simons Foundation, grant \#244930. }

\begin{abstract}
Let $G$ be a graph and let $I$ be the edge ideal of $G$. Our main results in this article provide lower bounds for the depth of the first three powers of $I$ in terms of the diameter of $G$. More precisely, we show that $\depth R/I^t \geq \left\lceil{\frac{d-4t+5}{3}} \right\rceil +p-1$, where $d$ is the diameter of $G$, $p$ is the number of connected components of $G$ and $1 \leq t \leq 3$. For general powers of edge ideals we show that $\depth R/I^t \geq p-t$. As an application of our results we obtain the corresponding lower bounds for the Stanley depth of the first three powers of $I$. 

\end{abstract}

\maketitle

\section{Introduction}
Let $R$ be either a Noetherian local ring or a standard graded $k$-algebra, where $k$ is a field. Let $I$ be an ideal of $R$ and when $R$ is graded assume that $I$ is a graded ideal. Let $d=\dim R$. A classical result by Burch \cite{Bur}, which was improved by Broadmann \cite{Bro}, states that $$\underset{t \rightarrow \infty}{\rm lim} \depth R/I^{t} \leq d-\ell(I),$$ where $\ell(I)$ is the analytic spread of $I$. Eisenbud and Huneke showed that the equality holds if the associated graded ring $\assg{I}= \bigoplus \limits_{i=0}^{\infty} I^i/I^{i+1}$ of $I$ is Cohen--Macaulay \cite{EH}. Therefore the limiting behavior of the depth is well understood. However the initial behavior of the depth of powers is still mysterious. Thus it is natural to investigate lower bounds for $\depth R/I^t$.

In the case of monomial ideals, lower bounds for the depth of the first power,  $\depth R/I$, have been studied extensively \cite{FHV, GV, K}. Herzog and Hibi  determined that $\depth R/I^{t}$ is a non--increasing function if all the powers of  $I$  have a linear resolution \cite{HH}. They also obtained lower bounds for $\depth R/I^t$ if all the powers of $I$ have linear quotients, a condition that implies that all the powers of $I$ have linear resolutions \cite{HH}. In particular, they showed that  all edge ideals associated to a finite graph whose complementary graph is chordal have linear quotients. Also, if $I$ is a square-free Veronese ideal (which includes the class of complete graphs) then all powers of $I$ have linear quotients. However, in general edge ideals and their powers do not  have linear resolutions. Even for monomial ideals the depth function can behave quite wildly, see \cite{BHH}. For square-free monomial ideals it is known that $\depth R/I^{t}$ will not necessarily be a non--increasing function, see \cite[Theorem~13]{KSS}, but the question is still open for edge ideals of graphs.

Another motivation for studying lower bounds for $\depth R/I^t$ is the fact that these lower bounds provide upper bounds for ${\rm projdim}_{R} R/I^t$, the projective dimension of $R/I^t$. When $I$ is the edge ideal of a graph then an upper bound for the projective dimension of a graph's edge ideal provides a lower bound for the first non-zero homology group of the graph's independence complex \cite[Observation~1.2]{DS}. Moreover, when $I$  is square-free monomial, its cohomological dimension and projective dimension are equal, \cite[Theorem~0.2]{EMS} or \cite[Corollary~4.2]{SW}. 
Many researchers have studied the question of finding upper bounds for the projective dimension of $R/I$ and upper bounds for the cohomological dimension, see for example \cite{EngNor, Fal, Hart, Lyu, Ogus, PS}.

We now describe our setup.  Let $V=\{x_1, \ldots, x_n\}$ be a set of $n$ vertices and let $G$ be a simple graph (no multiple edges, no loops) on $V$. Let $I$ be the edge ideal of $G$ in the ring  $R=k[x_1, \ldots, x_n]$, where $k$ is a field. By $\depth R/I^{t}$ we mean the maximum length of an $R/I^{t}$-regular sequence in $\m=(x_1, \ldots, x_n)$. When $I$ is the edge ideal of a bipartite graph then $\depth R/I^t \geq 1$, since $\m \not\in {\rm Ass}(R/I^t)$, by \cite[Theorem 5.9]{SVV}. In a recent article, Morey gives lower bounds for the depths of all powers when $I$ is the edge ideal of a forest, \cite{Mor}. 
We focus our interest on studying lower bounds for the depths of powers of edge ideals of graphs without any restrictions on the shape of the graph.

 The article is organized as follows. In Section~\ref{background} we give the necessary definitions and relevant background. In Sections~\ref{results} and~\ref{higher} we 
 establish the main results of this article. More precisely, we prove in Theorem~\ref{graph} that when $I$ is the edge ideal of a graph then $\depth R/I \geq
\left\lceil{\frac{d+1}{3}} \right\rceil $, where $d$ is the diameter of the graph. One can improve this bound by considering the diameters of each connected component of the graph. We show in  Corollary~\ref{sum connected comp} that when $G$ has $p$ connected components then $\depth R/I \geq \sum \limits_{i=1}^{p}\left\lceil{\frac{d_i+1}{3}}\right\rceil $, where $d_i$ is the diameter of the $i$-th connected component of $G$. 

We develop a series of lemmas that leads us to prove lower bounds for the second and third powers of the edge ideal of a graph. We first prove in Proposition~\ref{depth conn. comp} that $\depth R/I^t \geq p-t$ for any $t$, then in Theorems~\ref{graphSquare} and ~\ref{graphCube} we show that $\depth R/I^2 \geq \left\lceil{\frac{d-3}{3}}\right\rceil +p-1$ and $\depth R/I^3 \geq \left\lceil{\frac{d-7}{3}}\right\rceil +p-1$, where $I$ is the edge ideal of a graph $G$, $d$ is the diameter of $G$ and $p$ is the number of connected components of $G$. It is worth noting here that in order to establish the bounds for the second and third powers we need to deal with the depth of the edge ideal of a graph that potentially has loops. 
We provide a lower bound on the depth of the edge ideal of a graph with loops based on knowledge of the position of the loops. More precisely, we prove in  Proposition~\ref{LoopsLemma} that when $I$ is the edge ideal of a graph with loops and $\ell$ is an integer such that there exists a vertex $u$ with $d(u,x) \geq \ell$ for all vertices $x$ for which there is a loop on $x$, then $\depth R/I \geq \left\lceil{\frac{\ell-1}{3}}\right\rceil$. This result for the depth of the edge ideal of a graph with loops is of independent interest.

We conclude the article by using \cite[Proposition~2.6]{BKU} or \cite[Lemma~2.2]{Ra} in place of the Depth Lemma, to extend  our results  to provide lower bounds on the Stanley depth of the powers of $I$. In particular, in Theorem~\ref{Stanley} we show that $\sdepth {R/I^t} \geq \left\lceil {\frac{d-4t+5}{3}} \right\rceil +p-1$ for $1 \leq t \leq 3$, where $\sdepth$ denotes the Stanley depth.

\section{Background} \label{background}

Let $V=\{x_1, \ldots, x_n\}$ be a set of $n$ vertices and let $G$ be a graph on $V=V(G)$. Let $E=E(G)$ denote the set of edges of $G$. Unless otherwise stated we will assume that $G$ is a simple graph, that is, without loops and without multiple edges. Let $R=k[x_1, \ldots, x_n]$ be a polynomial ring, where $k$ is a field. Note that we will not distinguish between the vertices of a graph and the variables in the corresponding polynomial ring. The edge ideal $I(G)$ of a graph $G$ is defined to be the monomial ideal in the ring  $R$ generated by the monomials $x_ix_j$, where $\{x_i, x_j\} \in E$. Similarly, if $I$ is a square-free monomial ideal generated in degree two, $G(I)$ is the graph associated to $I$. That is, $\{x_i, x_j\} \in E(G(I))$ if and only if  $x_ix_j$ is a generator of $I$.

We now collect some useful definitions from graph theory. For algebraic definitions and background material, see \cite{mat} or \cite{Rafael}.

\begin{definition}{\rm 
Let $G$ be a graph, let $V=V(G)=\{x_1, \ldots, x_n\}$ and let $E=E(G)$. Then

\begin{enumerate} [(a)]
\item A {\it{path}} of length $r-1$ is a set of  $r$ distinct vertices $x_{i_1}, \ldots, x_{i_r}$ together with $r-1$ edges $x_{i_j}x_{i_{j+1}}$, where $x_{i_j} \in \{ x_1, \ldots, x_n\}$ and $1 \leq j \leq r-1$. 

\item The {\it distance} between two vertices $u$ and $v$ is  the length of the shortest path between $u$ and $v$ and is denoted $d(u,v)$. 

\item The {\it diameter} of a connected graph is $d(G)=\max \{d(u,v)\mid u, v \in V\}$. Therefore, if $d=d(G)$ then there exist vertices $u,v$ of $G$ with $d(u,v)=d$. In this case we say that a path of length $d$ with endpoints $u$ and $v$ {\it{realizes the diameter}} of $G$.  Although technically the diameter of a disconnected graph is infinite, we will find it useful to refer to the maximum of the diameters of the connected components of $G$ as the diameter of $G$ when $G$ is disconnected.

\item Let $u \in V$. The {\it{neighbor set}} of $u$ is the set $N(u)=\{v\in V(G) \mid \{u,v\} \in E\}$. When $N(u)=\emptyset$ then $u$ is called an {\it{isolated vertex}} and when the cardinality of $N(u)$ is one then $u$ is called a {\it{leaf}}.

\item A {\it{loop}} in a graph $G$ is an edge both of whose endpoints are equal, that is, an edge $\{x, x\} \in E$. A loop on $x$ corresponds to a generator $x^2$ in the edge ideal, so the edge ideal of a graph with loops is no longer square-free. Note that if loops are added to a graph, the distance between two vertices is unchanged.
\end{enumerate}
}
\end{definition}

When dealing with general graphs, it is helpful to consider a
construction that is commonly used to produce a spanning
tree. Although the spanning tree produced will not be used here, nonetheless the
construction yields a partition of the vertices that we will exploit.

\begin{notation} \label{XiSets}
{\rm Suppose $G$ is a connected graph and $u \in V(G)$. Define
$$X^i_G(u) =\{x\in V(G) \, | \, d(u,x)=i\}.$$ Note that $X^0_G(u)=\{u\}$ and that $i$ runs from $0$ to $d$, where $d=\mbox{\rm{max }}\{d(u,x) \, | \, x \in V(G)\}$. The sets $X^i_G(u)$ form a partition of $V(G)$. Once $u$ has been fixed, we will  omit $G$ and $u$ from the notation when they are clear from context.  We will frequently choose $u$ to be an endpoint of a path realizing the diameter, in which case $d$ will be the diameter of $G$. When $G$ is not connected, this construction can be applied to the connected component of $G$ containing $u$. 

When a vertex $u$ has been fixed in $G$, we will denote the connected component of $G$ that contains $u$ by ${}_{u}G$. Thus if $I$ is an edge ideal and $u$ has been fixed, then $d({}_{u}G(I))$ denotes the diameter of the connected component of $G(I)$ containing $u$.}
\end{notation}

There are two basic facts about these sets that will prove useful in the sequel. Fix $u$ and form $X^i=X^i_G(u)$.
First note that if $x \in X^i$ for $i \geq 1$, then $N(x) \cap X^{i-1}$ is nonempty since
there is a path from $u$ to $x$ of length precisely $i$ by the
definition of $X^i$. Also, if $u$ and $v$ are the endpoints of a path realizing the diameter, then $v \in X^d$ and if $y \in N(v)$, then $y$ is not a leaf. If $y$ were a leaf, $d(u,y)=d+1$, a contradiction.

The next lemma is well known, see for example \cite[Lemma~2.2]{Mor}.

\begin{lemma} \label{addVariable}
Let $I$ be an ideal in a polynomial ring $R$, let $x$ be an
indeterminate over $R$, and let $S=R[x]$. Then $\depth S/IS = \depth R/I +1$.
\end{lemma}

If $x_1$ is an isolated vertex of a graph $G$, define $R'=k[x_2, \ldots, x_n]$.  Notice that all generators of $I=I(G)$ lie in $R'$ and so by abuse of notation we can consider an ideal $I'=IR'$ in the ring $R'$ generated by the edges of the graph $G$. Then by Lemma~\ref{addVariable}, $\depth R/I = \depth R'/I' +1$. Thus we will assume graphs are initially free of isolated vertices and that all variables of $R$ divide at least one generator of $I$. 

Throughout the paper we will perform operations on ideals that correspond to the graph minors of contractions and deletions. A {\it deletion} minor is formed by removing a vertex $x$ from $G$ and deleting any edge of $G$ containing $x$. This corresponds to the ideal  $(I,x)$, or more precisely the quotient ring $R/(I,x)$. This process can result in isolated vertices, which will increase the depth of the quotient ring as in Lemma~\ref{addVariable}. To provide clarity we will count isolated vertices separately and will require connected components of a graph to have at least two vertices. A {\it contraction} minor of $G$ is formed by removing redundancies from the set $\{e\setminus \{x\} \mid e \in E(G)\}$ to obtain the edge set of the contraction. Note that if $y$ is a neighbor of $x$ in $G$, it becomes an isolated vertex of the contraction as any other edge containing $y$ was removed as a redundancy.
This corresponds to forming the ideal $(I:x)$. Note that $N(x) \subseteq (I:x)$ and so such an ideal may have variables as generators. However, if $K=(J,x_1)$ is a minimal generating set of an ideal, then $R/K \cong k[x_2, \ldots, x_n]/J$. Thus we will refer to $K$ as an edge ideal if $J$ is an edge ideal.

For clarity and ease of reference, we now state several previously known results.

\begin{lemma}\label{Generalization of HaMorey}
Let $I$ be a monomial ideal in a polynomial ring $R$ and
let $M$ be a monomial in $R$. If $y$ is a variable such that $y$ does
not divide $M$ and $K$ is the extension in $R$ of the image of $I$ in $R/y$, then $((I:M),y)=((K:M),y)$.
\end{lemma}

\proof
See the proof of  \cite[Theorem 3.5]{HaMorey}.
\QED

\begin{lemma}\cite[Lemma~2.10]{Mor}\label{leaf}
Suppose $G$ is a graph, $I=I(G)$, $x$ is a leaf of $G$, and $y$ is the
unique neighbor of $x$. Then $(I^t:xy)=I^{t-1}$ for any $t \geq 2$.
\end{lemma}

We conclude this section with an extension of the preceding lemma 
that will allow us to use any edge of the graph.

\begin{lemma}\label{edge}
Let $G$ be a graph, $I=I(G)$ and $\{x, y\}\in E(G)$. Then $(I^2:xy)=(I,E)$, where $E = \, <x_iy_j \, | \, x_i\in N(x)\, {\mbox{\rm and}}\, 
y_j\in N(y)>$. More generally, if $x_1\cdots x_{2t}\in I^{t}$, then $(I^{t+1}:x_1\cdots x_{2t})=(I,E)$, where $E$ is the ideal generated by all degree two monomials $y_1y_2$ supported on $\bigcup \limits_{i=1}^{2t} N(x_i)$ satisfying $y_1y_2x_1 \cdots x_{2t}\in I^{t+1}$.
\end{lemma}

\proof
Suppose first that $a$ is a minimal generator of $(I,E)$. If $a \in
I$, then $a \in (I^{2}:xy)$ since $xy \in I$. Else $a=x_iy_j\in E$ and
$axy=x_ixy_jy \in I^2$. Thus $(I,E) \subseteq (I^2:xy)$.

Conversely, suppose $b\in (I^2:xy)$ but $b \not\in I$. Since
$(I^2:xy)$ is a monomial
ideal, we may assume that $b$ is a monomial. Then $bxy \in I^2$, so
$bxy=e_1e_2h$, where $e_i$ are degree two monomials corresponding to
edges of $G$. Since $b \not\in I$, $e_i$ does not divide $b$ for
$i=1,2$, and so
without loss of generality, $x$ divides $e_1$ and $y$ divides
$e_2$. Thus $e_1=xx_i$ and $e_2=yy_j$ for some $x_i\in N(x)$ and $y_j
\in N(y)$. Thus $x_iy_j$ divides $b$ and so $b \in E \subset (I,E)$. 

The proof of the generalized statement follows the same outline. Note that the $x_i$ need not all be distinct.
\QED

Note that the ideal $(I,E)$ in Lemma~\ref{edge} is no longer guaranteed to be square-free. If $z\in
N(x) \cap N(y)$, then $z^2 \in E$. However, 
$(I,E)$ is still a monomial ideal, and if $z^2$ and $w^2$ are both
generators of $E$,
then $zw \in (I,E)$. This follows easily since $z\in N(x)$ and $w \in
N(y)$.

\section{The first power} \label{results}

As a first step toward determining the depths of $R/I^t$ for arbitrary
graphs, a lower bound, similar to the one given in \cite{Mor} for trees, is needed for $\depth R/I$. This lower bound is generally
far from sharp, however it is of a form that generalizes to higher powers. Alternate bounds for this depth, or equivalently for the projective dimension of $R/I$, exist in the literature, \cite{DHS,DS,DocEng,HH}. However the focus here is on providing a bound that will serve as the basis for bounds on the depths of higher powers, using techniques that will extend to higher powers. We first present the main result of this section. An alternate proof has been communicated to us by Russ Woodroofe.

\begin{theorem}\label{graph}
Let $G$ be a connected graph and let $I=I(G)$. If there exist $u, v \in V(G)$ with $d(u,v)=d$, then $\depth R/I \geq
\left\lceil{\frac{d+1}{3}} \right\rceil $.
\end{theorem}

\proof
We proceed by induction on $n$, the number of vertices. Notice that for any fixed $d$, we have that  $n\geq d+1$.
Since $\m \not\in \Ass(R/I)$, then $\depth R/I \geq 1$. Note that if
$d\leq 2$, then $\left\lceil{\frac{d+1}{3}}\right\rceil=1$ and so the result
holds. If $n=d+1$, the graph is a path and thus the result holds by \cite[Lemma~2.8]{Mor}. Hence we may assume $n-1>d \geq 3$. 

Let $X^{i}=X^i_G(u)$ be as in Notation~\ref{XiSets} 
and let $w \in N(v) \cap X^{d-1}$. 
Consider first $(I:w)=(J,N(w))$,  where $J$ is the ideal
corresponding to the minor $G'$ of $G$ formed by deleting the variables in
$N(w)$. 
Since $d\geq 3$ then $X_{G'}^{d-3}(u) \neq \emptyset$. Let $z\in X_{G'}^{d-3}(u)$ and 
notice that $d(u,z)=d-3$. Moreover, $w$ does not divide any generator of
$(J,N(w))$. Thus $(J,N(w)) \subset R'[N(w)]$, where $R'$ is the polynomial
ring formed by deleting $w\cup N(w)$. Then we have
\begin{eqnarray*}
\depth R/(I:w) &=& \depth R'[w,N(w)]/(J,N(w))\\
&=&\depth R'[w]/J =\depth R'/J+1 \\
& \geq & \left\lceil{\frac{d-3+1}{3}}\right\rceil +1 =
\left\lceil{\frac{d+1}{3}}\right\rceil
\end{eqnarray*}
by induction on $n$.

Next we consider $(I,w)=(K,w)$, where $K$ is the ideal of the minor $G''$
of $G$ formed
by deleting $w$. If $G''$ is connected, then $d(u,v)=d$ in $G''$ and therefore
$\depth R/(I,w)=\depth R/(K,w) \geq \left\lceil{\frac{d+1}{3}}\right\rceil$ by
induction on $n$. If $G''$ is not connected, then there is a vertex $z \in {}_{u}G''$ with $d(u,z) \geq d-2$ 
and $v \not\in {}_{u}G''$.
If $v$ is an
isolated vertex, then by Lemma~\ref{addVariable} we obtain $\depth R/(K,w) \geq
\left\lceil{\frac{d-2+1}{3}}\right\rceil +1 \geq
\left\lceil{\frac{d+1}{3}}\right\rceil$.
Otherwise, $v$ is in a connected
component of $G''$ that has depth at least one, so by \cite[Lemma~6.2.7]{Rafael}, we have $$\depth R/(K,w) \geq
\left\lceil{\frac{d-2+1}{3}}\right\rceil +1 \geq
\left\lceil{\frac{d+1}{3}}\right\rceil.$$ In either case, $\depth R/(I,w) \geq
\left\lceil{\frac{d+1}{3}}\right\rceil$.

Applying the Depth Lemma~\cite[Proposition 1.2.9]{BH} to the short exact
sequence
$$0 \rar R/(I:w) \rar R/I \rar R/(I,w) \rar 0$$
yields $\depth R/I \geq \left\lceil{\frac{d+1}{3}}\right\rceil$, as desired.
\QED

By selecting a pair of vertices $u$ and $v$ whose distance is maximal, we immediately obtain the following corollary.

\begin{corollary}\label{graphDiameter}
Let $G$ be a connected graph of diameter $d\geq 1$ and let $I=I(G)$. Then $\depth R/I \geq
\left\lceil{\frac{d+1}{3}} \right\rceil $.
\end{corollary}

As an immediate corollary we extend Theorem~\ref{graph} to graphs that are not necessarily connected.

\begin{corollary}\label{sum connected comp}
Let $G$ be a graph with $p$ connected components, $I=I(G)$, and let $d_i$ be the diameter of the $i$th connected component. 
Then $\depth R/I\geq \sum \limits_{i=1}^{p}\left\lceil{\frac{d_i+1}{3}}\right\rceil $. In particular, $\depth R/I \geq \left\lceil{\frac{d+1}{3}}\right\rceil +p-1$.
\end{corollary}

\proof
This follows directly from Theorem~\ref{graph} and \cite[Lemma~6.2.7]{Rafael}. \QED

The next corollary is an interesting result that follows from the proof of Theorem~\ref{graph}. Although the result could be used to prove the theorem above, it is difficult to obtain independently. However, it can be useful in bounding the depths of higher powers.

\begin{corollary} \label{I:w}
Let $G$ be a graph, let $I=I(G)$, and fix $u \in V(G)$. Let $w\in  X^{\ell}=X^{\ell}_G(u)$ for some $0\leq \ell$. Then $\depth R/(I:w) \geq \left\lceil{\frac{\ell+2}{3}}\right\rceil$.

\end{corollary}

\proof  Let $w\in X^\ell$. Notice that $(I:w)=(J, N(w))$, where $J$ is the ideal corresponding to the minor $G'$ of $G$ formed by deleting the variables in $N(w)$. Let $R'$ be the polynomial ring formed by deleting $w$ and the variables in $N(w)$. 
As before we have
\begin{eqnarray*}
\depth R/(I:w) &=& \depth R'[w,N(w)]/(J,N(w))\\
&=&\depth R'[w]/J =\depth R'/J+1. 
\end{eqnarray*}
If $\ell <2$ the result holds since  $\depth R/(I:w) \geq 1$. Hence we may assume that $\ell \geq 2$. Since the diameter of $G'$ is at least $\ell-2$, applying Theorem~\ref{graph} yields 
$$\depth R'/J+1\geq \left\lceil{\frac{\ell-2+1}{3}}\right\rceil +1 =\left\lceil{\frac{\ell+2}{3}}\right\rceil $$
and the result follows.
\QED

We conclude this section with an extension of Theorem~\ref{graph} that gives a bound for the depth of the first power of the edge ideal of a graph with loops. This result is of independent interest.

\begin{proposition}\label{LoopsLemma}
Let  $G$ be a connected graph with loops and let $I=I(G)$. If there exists $u\in V(G)$ with $d(u,x) \geq \ell$ for all $x$ such that $\{x,x\} \in E(G)$, then $\depth R/I \geq \left\lceil {\frac{\ell-1}{3}} \right\rceil$.
\end{proposition}

\proof
Notice that if $\ell <2$ the result is trivial. Thus we assume $\ell \geq 2$. We induct on the number of loops.  Let $x$ be a variable corresponding to a vertex with a loop. Notice that $(I:x)=(I,N(x))=(J,N(x))$, where $J$ is the minor formed by deleting all vertices in $N(x)$. 
Since $x \in N(x)$, the number of loops of $G(J)$ is less than the number of loops of $G$.  Notice that since all deleted vertices are at least distance $\ell -1$ from $u$,  $d({}_{u}G(J)) \geq \ell -2$ and $d(u,z) \geq \ell$ for all loops $z$. If ${}_{u}G(J)$ has no loops, then  $\depth R/(I:x) \geq \left\lceil {\frac{\ell-1}{3}} \right\rceil$, by  Theorem~\ref{graph} since $\left\lceil {\frac{d(G(J))+1}{3}} \right\rceil \geq \left\lceil {\frac{\ell-2+1}{3}} \right\rceil$.  
If ${}_{u}G(J)$ contains a loop $y$, then 
$d(u,y) \geq \ell$ and hence $\depth R/(I:x) \geq \left\lceil {\frac{\ell-1}{3}} \right\rceil$, by induction.

Now consider $(I,x)=(K,x)$, where $K$ is the minor of $I$ formed by deleting $x$. 
Then $d({}_{u}G(K)) \geq \ell-1$ and $G(K)$ has fewer loops than $G$, so $\depth R/(I,x) \geq \left\lceil {\frac{\ell-1}{3}} \right\rceil$, by either Theorem~\ref{graph} or induction as above.

Applying the Depth Lemma~\cite[Proposition 1.2.9]{BH} to the short exact sequence 
$$0 \rar R/(I:x) \rar R/I \rar R/(I,x) \rar 0.$$
completes the proof. \QED

\section{Depths of Higher Powers of Edge Ideals}\label{higher}

Our main results in this section focus primarily on $I^2$ and $I^3$. Selected results are stated for all powers since our methods can extend to higher powers, particularly when one has some control over the structure of the underlying graph. 
The central idea of the proofs will be to apply the Depth Lemma~\cite[Proposition 1.2.9]{BH} to families of short exact sequences. We begin the section by introducing some notation.

We will frequently use deletion minors in the proofs, and often the minors will be formed using a collection of vertices.  Let $G$ be a graph and let $I=I(G)$. For $a \in V(G)$ we let $I_a$ represent the edge ideal of the minor of $G$ formed by deleting $a$. We will refer to $I_a$ as a minor of $I$. Given a collection of vertices $y_1, \ldots , y_s$, define $I_0=I$ and for $1 \leq i \leq s$ define $I_i$ to be the minor of $I$ formed by deleting $y_1, \ldots , y_i$. Define $R_i$ to be the corresponding polynomial ring, namely $R_i=R/(y_1, \ldots, y_i)$. 

Recall  that an {\it{ induced graph}} on a subset  $\{x_1, \ldots, x_r\}$ of vertices of a graph $G$ is a graph $G'$ with $V(G')=\{x_1, \ldots, x_r\}$ and $E(G')=\{ \{x_{i},x_{j}\} \in E(G) \mid x_i, x_j \in V(G')\}$.

\begin{lemma}\label{OrderNeighbors}
Let $G$ be a graph, $V=V(G)$ and $I=I(G)$. Let $x_1,\ldots,  x_r \in V$ be such that the induced graph on $x_1, \ldots, x_r$ is connected and fix a vertex $u$ in the connected component of $G$ containing $x_1, \ldots , x_r$. Let $\{y_1, \ldots ,y_s\} \subset \bigcup \limits_{i=1}^{r} N(x_i) \setminus \{x_1, \ldots ,x_r\}$. Then there exists an ordering of the vertices $y_1, \ldots ,y_s$ such that for all $i <s$, $x_1, \ldots , x_r \in {}_{u}G(I_i)$, 
where $I_i$ is obtained by deleting $y_1, \ldots, y_i$. 
\end{lemma}

\proof Using the fixed vertex $u$, form $X^i=X^i_G(u)$. Notice that $u$ may be one of the $y_i$. Since $x_1, \ldots ,x_r \in{}_{u}G$, then for each $i$, $x_i \in X^t$ for some $t$. Let $k$ be the least positive integer for which $x_i \in X^k$ for some $i$. Fix $x_q\in X^k$. Then there is a path from $u$ to $x_q$ containing precisely one vertex in $X^j$ for each $j \leq k$. Since for every $i$, $y_i \in N(x_{\ell})$ for some $\ell$, then $y_i \in \bigcup \limits_{j=k-1}^d X^j$ for all $i$. Thus at most one $y_i$ lies on the chosen path.  We may reorder the variables so that $y_s$ is this vertex (if any). Then for all $i<s$,   there is a path in $I_i$ from $u$ to $x_q$ and there is a path from $x_q$ to $x_i$ for all other $i$, since the induced graph on $x_1, \ldots, x_r$ is connected. \QED

Once we have ordered a collection of neighboring vertices as in Lemma~\ref{OrderNeighbors}, deleting the vertices in order will result in a series of graphs for which $u$ and $x_1 , \ldots , x_r$ are in the same connected component, followed by a graph for which $u$ and $x_i$ might be disconnected. When $r=1$ and $\{y_1, \ldots, y_s\} = N(x_1)$, deleting all vertices except $y_s$ will result in a graph for which $x_1$ is a leaf. The next lemma formalizes how this can be used to estimate depths. Although it will generally be used when $M=x_1$ is a single vertex or $M=x_1\cdots x_r$ is the product of connected vertices and $\{y_1, \ldots, y_s\}=N(x_r)\setminus \{x_1, \ldots ,x_{r-1}\}$, the result holds in the more general situation described here.

\begin{lemma}\label{ExhaustTheNeighbors}
Let $R$ be a polynomial ring over a field, $I$ an ideal, and let $M$ be a monomial in $R$. Let $\{y_1, \ldots, y_s\}$ be variables such that for all $i$, $y_i$ does not divide $M$. Let $a, b$ be two nonnegative integers. If $\depth R_{i-1}/(I_{i-1}^t:My_{i}) \geq a$ for all $i\geq 1$ and $\depth R_s/(I_s^t:M) \geq b$, then $\depth R_i/(I_i^t:M) \geq \min \{a,b\}$ for each $i\geq 0$. In particular, $$\depth R/(I^t:M) \geq \min \{a, b\}.$$
\end{lemma}

\proof
Consider the family of short exact sequences 
\begin{eqnarray*}
0 & \rar & R/(I^t:My_1) \rar R/(I^t:M) \rar R/((I^t:M),y_1) \rar 0 \\
0 & \rar & R_1/(I_1^t:My_2) \rar R_1/(I_1^t:M) \rar R_1/((I_1^t:M),y_2) \rar 0 \\
0 & \rar & R_2/(I_2^t:My_3) \rar R_2/(I_2^t:M) \rar R_2/((I_2^t:M),y_3) \rar 0 \\
& \vdots & \\
0 & \rar & R_{s-1}/(I_{s-1}^t:My_s) \rar R_{s-1}/(I_{s-1}^t:M) \rar R_{s-1}/((I_{s-1}^t:M),y_s) \rar 0.
\end{eqnarray*}
Notice that by Lemma~\ref{Generalization of HaMorey} the right hand term of sequence $i$ is isomorphic to $R_i/(I_i^t:M)$, which is the center term of sequence $i+1$. Now $\depth R_i/(I_i^t:My_i) \geq a$ by hypothesis and 
$R_{s-1}/((I_{s-1}^t:M),y_s) \cong R_s/(I_s^t :M)$, so by hypothesis, $\depth R_{s-1}/((I_{s-1}^t:M),y_s) \geq b$. By applying the Depth Lemma~\cite[Proposition 1.2.9]{BH} repeatedly starting with the final sequence and working our way up we see that $\depth R_i/(I_i^t:M) \geq \min \{a,b\}$ for each $i$ from $i=s-1$ to $i=0$. Since $\depth R_s/(I_s^t:M) \geq b$, the result holds for all $i$.
\QED

We now give a first estimate on the depth of any power of an edge ideal in terms of the number of connected components of the graph. Recall that we have defined connected components to have at least two vertices. In Corollary~\ref{sum connected comp} we were able to achieve a better bound for the first power and later in this section we will improve this bound for the second and third powers; however, the advantage of considering this bound is that it is a bound for all the powers even though it might not be sharp. 

\begin{proposition}\label{depth conn. comp}
Let $G$ be a graph with $p$ connected components and let $I=I(G)$. Then for every $t\geq 1$ $$\depth R/I^t \geq p-t.$$
\end{proposition}

\proof
We prove this by induction on $p$, the case of $p=1$ being clear. Suppose that $p\geq 2$ and $I=(J,K)$, where $J \subset A=k[x_1, \ldots, x_r]$ is the edge ideal of the graph consisting of all but one of the connected components of $G$  and $K \subset B=k[x_{r+1}, \ldots, x_n]$ is the edge ideal of the remaining connected component of $G$. Then $\depth A/J \geq p-1$ and $\depth B/K \geq 1$ by Corollary~\ref{sum connected comp}. By induction on $p$ we have $\depth A/J^s \geq p-1-s$ for all $s\geq 1$. In particular, $\depth A/J^{t-i} \geq p-1-(t-i)=p-t+i-1$ for $1 \leq i \leq t-2$ and for all $1\leq j \leq t$, $\depth A/J^{t-j+1} \geq p-1-(t-j+1)=p-t+j-2$. Then by \cite[Theorem~2.4]{HTT} we have 
\begin{eqnarray*}
\depth R/I^t &\geq& \min_{\begin{subarray}{l}{i \in [1,t-1]}\\{j\in [1,t]}\end{subarray}} \{ \depth A/J^{t-i}+\depth B/K^i+1, \depth A/J^{t-j+1}+\depth B/K^{j}\}\\
&=&\min_{\begin{subarray}{l}{i \in [1,t-2]}\\{j\in [2,t]}\end{subarray}}\{ \depth A/J^{t-i}+\depth B/K^i+1, \depth A/J+\depth B/K^{t-1}+1, \\
&&\depth A/J^{t}+\depth B/K,\depth A/J^{t-j+1}+\depth B/K^{j}\}  \\
&=&  \min_{\begin{subarray}{l}{i \in [1,t-2]}\\{j\in [2,t]}\end{subarray}} \{p-t+i-1+0+1, p-1+0+1, p-t-1+1, p-t+j-2+0 \}\\
&=&  \min _{\begin{subarray}{l}{i \in [1,t-2]}\\{j\in [2,t]}\end{subarray}} \{p-t+i, p, p-t, p-t+j-2 \}=p-t.
\end{eqnarray*}
\QED

The next theorem establishes a sharper lower bound for the  depth of the second power of an edge ideal.

\begin{theorem}\label{graphSquare}
Let $G$ be a graph with $p$ connected components, $I=I(G)$, and let $d=d(G)$ be the diameter of $G$.
Then
$$\depth R/I^2 \geq \left\lceil{\frac{d-3}{3}}\right\rceil+p-1.$$
\end{theorem}

\proof
We proceed by induction on $n$, the number of vertices in $G$. Suppose $n\leq 4$. Then $d \leq 3$ and $p \leq 2$ since the number of
connected components does not include isolated vertices. If $p=1$ the
bound is trivial. If $p=2$, for $n\leq 4$ the graph must be a forest
consisting of two disconnected edges and the result follows from \cite[Theorem~3.4]{Mor}. Note that in general, if $p \geq 2$, then
$\depth R/I^2 \geq 1$ by \cite[Lemma~2.1]{AJ}.

We may now assume that $n\geq 5$.
Let $u,v$ be the endpoints of a path that realizes the diameter and let $X^i=X^i_G(u)$. Let $w \in N(v)$ and let $\{y_1, \ldots, y_s\}=N(w)$ be ordered as in Lemma~\ref{OrderNeighbors} so that $d(u,w)$ is finite in $I_i$ for $i<s$. Recall that $I_0=I$. Then for each $1 \leq i \leq s$ we have $(I_{i-1}^2:wy_i)=(I_{i-1},E_{i-1})$, where $E_{i-1}$ is as in Lemma~\ref{edge}. Now $(I_{i-1},E_{i-1})$ is the edge ideal of a graph $G'$, possibly with loops, of diameter at 
least $d-1$ since $d(u,w) \geq d-1$ even with the additional edges. Thus if $(I_{i-1},E_{i-1})$ is square-free, $$\depth R_{i-1}/(I_{i-1},E_{i-1})\geq \lceil{\frac{d(G')+1}{3}}\rceil +p-1\geq \lceil{\frac{d-1+1}{3}}\rceil +p-1$$ by Corollary~\ref{sum connected comp}. If $(I_{i-1},E_{i-1})$ is not square-free,
then there exists $x\in V(G)$ such that  $x^2 \in E_{i-1}$. Now $x\in N(w)$, and so $d(u,x) \geq d-2$. Note that each connected component of $G((I_{i-1},E_{i-1}))$ other than ${}_{u}G((I_{i-1},E_{i-1}))$ will be square-free, and so have depth at least one. Thus combining \cite[Lemma~6.2.7]{Rafael} with Proposition~\ref{LoopsLemma} yields
$\depth R_{i-1}/(I_{i-1},E_{i-1}) \geq \lceil{\frac{d-2-1}{3}}\rceil+p-1$ for $i \leq s$. 

Now $w$ is isolated in $I_s$, so $(I_s^2:w)=I_s^2$ and $w$ is a free variable in $R_s/I_s^2$. 
Since $d({}_{u}G(I_s)) \geq d-3$, then by induction and  Lemma~\ref{addVariable} we have $$\depth R_s/(I_s^2:w)
\geq \lceil{\frac{d({}_{u}G(I_s))-3}{3}}\rceil +p-1+1\geq 
\lceil{\frac{d-3}{3}}\rceil +p-1.$$ Hence by Lemma~\ref{ExhaustTheNeighbors} we obtain $\depth R/(I^2:w) \geq
\lceil{\frac{d-3}{3}}\rceil +p-1$.

Finally, consider $(I^2,w)=(I_w^2,w)$. If $v \in  {}_{u}G(I_w)$, 
then $d({}_{u}G(I_w)) \geq d$ and $\depth R/(I^2,w) =\depth R_w/I_w^2 \geq
\lceil{\frac{d-3}{3}}\rceil+p-1$ by induction on $n$. 
Otherwise $d( {}_{u}G(I_w)) \geq d-2$ and $G(I_w)$ contains an additional connected component or an isolated vertex, so
\begin{eqnarray*}
\depth R/(I^2,w)&=& \depth R/(I_w^2,w) \\
&\geq & \left\lceil{\frac{d-2-3}{3}}\right\rceil+(p+1)-1 =
\left\lceil{\frac{d-2}{3}}\right\rceil+p-1.
\end{eqnarray*}

By applying the Depth Lemma~\cite[Proposition 1.2.9]{BH} to the following exact sequence 
$$0 \rar R/(I^2:w) \rar R/I^2 \rar R/(I^2,w)\rar 0$$ 
we see that $\depth R/I^2\geq
\left\lceil{\frac{d-3}{3}}\right\rceil+p-1$ as desired.
\QED

\begin{remark} \label{remarkdistance} {\rm
Notice that in the proof of Theorem~\ref{graphSquare} we required that $u$ and $v$ be endpoints of a path that realizes the diameter. This was done in order to obtain the best possible lower bound for the depth of $R/I^2$. However, one may take $u$ and $v$ to be endpoints of any path of length $\ell=d(u,v)$. Then continuing as in the proof of Theorem~\ref{graphSquare}, we would obtain that $\depth R/I^2\geq
\left\lceil{\frac{\ell-3}{3}}\right\rceil+p-1$. Although this is a weaker lower bound, it can be useful in a more general setting. }\end{remark}

As with the proof of Theorem~\ref{graph} the proof of Theorem~\ref{graphSquare} yields the following interesting corollary.

\begin{corollary} \label{Square Colon One}
Let $G$ be a graph 
and let $I=I(G)$. Fix $u \in V(G)$ 
and let $w \in X^{\ell}=X_G^{\ell}(u)$ for some $0 \leq \ell $. Then $\depth R/(I^2:w) \geq \left\lceil {\frac{\ell-2}{3}} \right\rceil$.
\end{corollary}

\proof First notice that when $\ell <2$ there is nothing to show. Hence we may assume that $\ell \geq 2$. 
Let $\{y_1, \ldots , y_s \}= N(w)$ be ordered as in Lemma~\ref{OrderNeighbors}. As in the proof of Theorem~\ref{graphSquare}, for each $1 \leq i \leq s$ we have  $(I_{i-1}^2:wy_i)=(I_{i-1},E_{i-1})$ as in Lemma~\ref{edge} and $(I_{i-1},E_{i-1})$ is the edge ideal of a graph of diameter at least $\ell$ since $d(u,w) =\ell$. Thus if $(I_{i-1},E_{i-1})$ is square-free, $\depth R_{i-1}/(I_{i-1},E_{i-1})\geq \lceil{\frac{\ell+1}{3}}\rceil +p-1$ by Corollary~\ref{sum connected comp}. If $x^2 \in E_{i-1}$, then $d(u,x) \geq \ell -1$ so combining \cite[Lemma~6.2.7]{Rafael} with Proposition~\ref{LoopsLemma} yields
$\depth R_{i-1}/(I_{i-1},E_{i-1}) \geq \lceil{\frac{\ell -2}{3}}\rceil$ for $i \leq s$. 

Now $w$ is isolated in $I_s$, and $d({}_{u}G(I_s)) \geq \ell -2$, so by Lemma~\ref{addVariable} and Theorem~\ref{graphSquare} 
$$
\depth R_s/(I_s^2:w) = \depth R_s/I_s^2+1 \geq 
\left\lceil{\frac{\ell -2-3}{3}}\right\rceil +1 = \left\lceil{\frac{\ell -2}{3}}\right\rceil.
$$
Hence by Lemma~\ref{ExhaustTheNeighbors} we have $\depth R/(I^2:w) \geq
\left\lceil{\frac{\ell -2}{3}}\right\rceil$. 
\QED

When exhausting the neighbors as in Lemma~\ref{ExhaustTheNeighbors}, we might end up with disconnected graphs. If the vertex $w$ is not in the connected component containing $u$, and thus is not in $X^i$ for any $i$, the bound above needs to be modified, but can still be found using only the diameter of ${}_{u}G(I)$.

\begin{lemma}\label{I^2:disconnected}
Let $G$ be a graph and let $I=I(G)$. Fix $u \in V(G)$ 
and let $w \in V(G)$ be such that $w\not\in  {}_{u}G$. Then  $\depth R/(I^2:w) \geq \left\lceil {\frac{\ell}{3}} \right\rceil$, where $\ell=d({}_{u}G)$.
\end{lemma}

\proof
Suppose $I=(J,K)$, where $K=I(_wG)$. Let $\{z_1, \ldots ,z_s\}$ be the neighbors of $w$ ordered as in Lemma~\ref{OrderNeighbors}. Note that $I_{i}=(J_i,K_i)$ and $J_i=J$ for all $i$. As in Lemma~\ref{edge}, we have $(I_{i-1}^2:wz_i)=(I_{i-1},E_{i-1})$, where all the edges in $E_{i-1}$ have endpoints in $V(_wG)$. Recall that $R_{i-1}$ is the polynomial ring corresponding to $I_{i-1}$ and let $R_{i-1}'$ be the polynomial ring with variables corresponding to $V(G(J_{i-1}))$. 
Then $\depth R_{i-1}/(I_{i-1},E_{i-1}) \geq \depth R'_{i-1}/J_{i-1} \geq\left\lceil {\frac{\ell +1}{3}} \right\rceil$, by \cite[Lemma~6.2.7]{Rafael} and Theorem~\ref{graph}. Finally, $w$ is an isolated vertex in $I_s$, so $(I_s^2:w)=I_s^2$ and $w$ is a free variable. Thus $\depth R_s/(I_s^2:w)\geq \depth R_s/I_s^2+1\geq 
 \left\lceil {\frac{\ell-3}{3}} \right\rceil+1 =  \left\lceil {\frac{\ell}{3}} \right\rceil$. The result then follows from Lemma~\ref{ExhaustTheNeighbors}.
\QED

 The lower bound for the depth of the first power of edge ideals that we obtained in Theorem~\ref{graph} is realized by edge ideals of paths, as was shown in \cite[Lemma~2.8]{Mor}. Therefore, one can not hope for any improvement of this bound for a general graph in terms of the invariants used. However, the lower bound for the depth of higher powers of edge ideals of paths given in \cite[Proposition~3.2]{Mor} is too high for general graphs. The next example shows that the bound we established in Theorem~\ref{graphSquare} is indeed attained, thus establishing that one can not improve this bound in terms of the invariants used.

\begin{example}\label{square sharp ex}
{\rm
Let $R=k[x_1, \ldots, x_{5}]$ and let $I$ be the edge ideal of the graph $G$ below

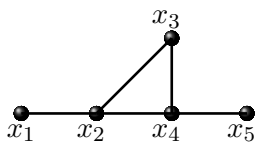
\begin{center}
\begin{tikzpicture}
\shade [shading=ball, ball color=black] (-1,0) circle (0.1) node [below ] {\textsf{\textbf{\textcolor{black}{\small{$x_1$}}}}};
\shade [shading=ball, ball color=black] (0,0) circle (0.1) node [below ] {\textsf{\textbf{\textcolor{black}{\small{$x_2$ }}}}};
\shade [shading=ball, ball color=black] (1,1) circle (0.1) node [above ] {\textsf{\textbf{\textcolor{black}{\small{$x_3$ }}}}};
\shade [shading=ball, ball color=black] (1,0) circle (0.1) node [below ] {\textsf{\textbf{\textcolor{black}{\small{$x_4$ }}}}};
\shade [shading=ball, ball color=black] (2,0) circle (0.1) node [below ] {\textsf{\textbf{\textcolor{black}{\small{$x_5$ }}}}};

\draw [line width=1pt  ] (-1,0)--(0,0);
\draw [line width=1pt  ] (0,0)--(1,1);
\draw [line width=1pt  ] (0,0)--(1,0);
\draw [line width=1pt  ] (1,0)--(1,1);
\draw [line width=1pt  ] (1,0)--(2,0);

\end{tikzpicture}

\end{center}

Then $d(G)=3$ and using Macaulay~2 \cite{M2} we have that $\depth R/I^2 =\lceil{\frac{d-3}{3}}\rceil=0$, which also follows from \cite[Theorem 3.3]{AJ}. Therefore, the bound in Theorem~\ref{graphSquare} is sharp.  

}
\end{example}

We now prove a series of lemmas that will allow us to establish a bound for the depth of the third power. Where possible, we give a general statement that holds for all powers. The first lemma is an extension of Lemma~\ref{edge}.

\begin{lemma}\label{CubeColonFour}
Let $G$ be a graph and let $I=I(G)$. Let $u, x_1, \ldots, x_{2t} \in V(G)$ for some $t \geq 1$ with $x_1\cdots x_{2t} \in I^t$. If for some $0 \leq \ell \leq d$ we have $x_i \in \bigcup \limits_{j=\ell}^d X^j$ for all $i$, where $X^j=X^j_G(u)$ 
then $\depth R/(I^{t+1}:x_1\cdots x_{2t}) \geq \left\lceil {\frac{\ell -2}{3}} \right\rceil$.
\end{lemma}

\proof
Notice that $(I^{t+1}:x_1\cdots x_{2t})=(I,E)$, where $E$ is the ideal generated by all degree two monomials $y_1y_2$ supported on $\bigcup \limits_{i=1}^{2t} N(x_i)$ satisfying $y_1y_2x_1\cdots x_{2t}\in I^{t+1}$ by Lemma~\ref{edge}. 
Let $G'$ be the graph, possibly with loops, associated to $(I,E)$.  Notice that $X^i_G(u)=X^i_{G'}(u)$ for $i \leq \ell -2$ since both endpoints of any generator of $E$ lie in $\bigcup \limits_{i=\ell -1}^d X^i_G$. This also implies that all loops of $G'$ are contained in $\bigcup \limits_{i=\ell -1}^d X^i_G$. So by Proposition~\ref{LoopsLemma} we have $\depth R/(I,E) \geq \left\lceil {\frac{\ell -2}{3}} \right\rceil$. 
\QED

The general outline of the following lemmas is to at each stage reduce by one the number of variables with which a colon ideal is formed. In general, this is accomplished using Lemma~\ref{ExhaustTheNeighbors}, however one must first deal with the situation where there are no neighbors to exhaust. This occurs when the graph is disconnected and one component consists of the induced graph on the variables used to form the colon ideal. In examining the later proofs in which the result is used, one sees that the goal is to create a path of vertices. The difficult case will be when the induced graph on the vertices of the path does not contain a leaf. Thus we assume in the next lemma that the graph contains a Hamiltonian cycle, that is, a cycle that passes through each vertex precisely once. To simplify notation, we will at times use ${\underline{x}}$ in place of $x_1, \ldots ,x_n$ when the number of variables used is clear.

\begin{lemma}\label{TriangleCase}
Let $G$ be a disconnected graph and let $I=I(G)$. Suppose $I=(J,K)$, where $J\subset k[x_1, \ldots ,x_n]$, $K \subset k[y_1, \ldots , y_{2t-1}]$, and $G(K)$ contains a Hamiltonian cycle. Then $\depth R/(I^{t+1}:y_1y_2\cdots y_{2t-1}) \geq \left\lceil {\frac{d(J)-3}{3}} \right\rceil$, where $d(J)=d(G(J))$ and $R=k[{\underline{x}},{\underline{y}}]$. 
\end{lemma}

\proof
Let $M=\prod_{i=1}^{2t-1}y_i$ and consider the family of short exact sequences

\begin{eqnarray*}
0  &\rar  &R/(I^{t+1}:My_1) \rar R/(I^{t+1}:M) \rar R/((I^{t+1}:M),y_1) \rar 0 \\
0  &\rar  &R/(((I^{t+1}:M),y_1):y_2) \rar R/((I^{t+1}:M),y_1) \rar R/((I^{t+1}:M),y_1,y_2) \rar 0 \\
 & \vdots & \\
0  &\rar & R/(((I^{t+1}:M),{\underline{y}}'):y_{2t-1}) \rar R/((I^{t+1}:M),{\underline{y}}') \rar R/((I^{t+1}:M),{\underline{y}}) \rar 0,
\end{eqnarray*}
where ${\underline{y}}'=\{y_1, \ldots, y_{2t-2}\}$.

We first handle the left hand term of each sequence by showing that for each $0 \leq i \leq 2t-2$, $(((I^{t+1}:M),y_1,\ldots , y_{i}):y_{i+1})=(J,K_i)$ for some $K_i$ an ideal of $k[y_1, \ldots , y_{2t-1}]$. Here we define $(((I^{t+1}:M),y_0):y_{1})=(I^{t+1}:My_1)$ since there is no element $y_0$.  By Lemma~\ref{Generalization of HaMorey} and some straight forward computations we have
$$(((I^{t+1}:M),y_1, \ldots , y_i):y_{i+1})=((I^{t+1}:My_{i+1}),y_1, \ldots, y_i).$$ Now since $M$ is a product of $2t-1$ variables that form a cycle and $y_{i+1}$ is an element of the cycle, $My_{i+1} \in I^t$ for each $i$. Thus by Lemma~\ref{edge}, $((I^{t+1}:My_{i+1}),y_1, \ldots , y_i)=(I,E_i, y_1, \ldots , y_i)$, where $E_i$ is the ideal generated by all degree two monomials $y_{i_1}y_{i_2}$ supported on $\bigcup \limits_{i=1}^{2t-1} N(y_i)$ satisfying $y_{i_1}y_{i_2}My_{i+1}\in I^{t+1}$. Now $M=\prod_{i=1}^{2t-1}y_i$ and $N(y_i) \subset k[y_1, \ldots , y_{2t-1}]$ for each $i$, so $(I,E_i,y_1, \ldots , y_i)=(J,K_i)$ where $K_i \subset k[y_1, \ldots , y_{2t-1}]$. Thus by \cite[Lemma~6.2.7]{Rafael} and Theorem~\ref{graph},
\begin{eqnarray*}
\depth R/((I^{t+1}:M),y_1,\ldots , y_{i}):y_{i+1}) &=&\depth k[{\underline{x}}]/J + \depth k[{\underline{y}}]/K_i \\
 & \geq & \depth k[{\underline{x}}]/J \geq \left\lceil {\frac{d(J)-3}{3}} \right\rceil.
\end{eqnarray*}

Now we claim $((I^{t+1}:M),y_1,\ldots , y_{2t-1})=(J^2,y_1 \ldots , y_{2t-1})$. Since $M\in I^{t-1}$, one inclusion is clear. Suppose $aM\in I^{t+1}$ for some monomial $a \not\in J^2$. Then $aM= e_1\cdots e_{2t+1}h$ for some monomial $h$ and some edges $e_i$. Note that if $e_i \in k[{\underline{x}}]$, then $e_i \mid a$ since $M \in k[{\underline{y}}]$. Since $a\not\in J^2$, at most one edge $e_i$ is in $ k[{\underline{x}}]$. Thus the ${\underline{y}}$-degree of  $e_1\cdots e_{2t+1}h$ is at least $2t$ but the degree of $M$ is $2t-1$ and hence $y_i \mid a$ for some $i$. Thus $(I^{t+1}:M) \subseteq (J^2,y_1 \ldots , y_{2t-1})$ and the second inclusion follows.  Thus 
\begin{eqnarray*}
\depth R/((I^{t+1}:M),y_1,\ldots , y_{2t-1}) & = & \depth R/(J^2,y_1 \ldots , y_{2t-1}) \\
 & = & \depth k[{\underline{x}}]/J^2 \geq \left\lceil {\frac{d(J)-3}{3}} \right\rceil
\end{eqnarray*} 
by Theorem~\ref{graphSquare}.
The result now follows from repeated applications of the Depth Lemma~\cite[Proposition 1.2.9]{BH}.
\QED

We now return to our computations concerning the depths of various ideals involving the third power of an edge ideal. 

\begin{lemma}\label{CubeColonThree}
Let $G$ be a graph and let $I=I(G)$. Let $u, x_1,x_2,x_3 \in V(G)$ and suppose that that $x_1,x_3 \in N(x_2)$ and $x_1,x_2,x_3 \in \bigcup \limits_{i=\ell}^dX^i$, where $X^i=X^i_G(u)$ for some $0\leq \ell \leq d$. Then $\depth R/(I^3:x_1x_2x_3)\geq \left\lceil{\frac{\ell -5}{3}}\right\rceil$.
\end{lemma}

\proof
We may assume $\ell \geq 6$ since otherwise the bound is trivial. First suppose $x_3$ is a leaf. Then $(I^3:x_1x_2x_3)=(I^2:x_1)$ and by Corollary~\ref{Square Colon One} we have $$\depth R/(I^3:x_1x_2x_3)=\depth R/(I^2:x_1) \geq \left\lceil{\frac{\ell-2}{3}}\right\rceil.$$

Suppose $x_3$ is not a leaf. We consider two cases. If $x_1x_3$ is a generator of $I$, let $\{z_1, \ldots, z_s\}=N(x_1)\cup N(x_2) \cup N(x_3) \setminus \{x_1,x_2,x_3\}$. If $x_1x_3$ is not a generator of $I$, let $\{z_1, \ldots, z_s\}=N(x_3)\setminus \{x_2\}$. In either case, order the vertices $z_1, \ldots, z_s$ as in Lemma~\ref{OrderNeighbors}.  Then by considering ${}_{u}G(I_{i-1})$, we have  $\depth R_{i-1}/(I_{i-1}^3:x_1x_2x_3z_i) \geq  \left\lceil{\frac{\ell -3}{3}}\right\rceil$ by Lemma~\ref{CubeColonFour} since $z_i \in \bigcup \limits_{i=\ell -1}^dX^i$. If $x_1x_3 \in I$, then $x_1,x_2,x_3$ forms a Hamiltonian cycle of a component that is disconnected from ${}_{u}G(I_s)$, so by Lemma~\ref{TriangleCase} we have that 
$\depth R_s/(I_s^3:x_1x_2x_3) \geq \left\lceil{\frac{d(I_s)-3}{3}}\right\rceil \geq \left\lceil{\frac{\ell -5}{3}}\right\rceil$, since $d({}_{u}G(I_s)) \geq \ell-2$.
When $x_1x_3 \not\in I$ then $x_3$ is a leaf in $I_s$, so as above, $(I_s^3:x_1x_2x_3)=(I_s^2:x_1)$. If $I_s$ is disconnected, then $d({}_{u}G(I_s))\geq \ell -2$. Thus by Lemma~\ref{I^2:disconnected}, or Corollary~\ref{Square Colon One} when ${}_{u}G(I_{s})$ is connected, we obtain
$\depth R_s/(I_s^3:x_1x_2x_3)\geq \left\lceil{\frac{\ell -2}{3}}\right\rceil$.
In either case, applying Lemma~\ref{ExhaustTheNeighbors} yields $\depth R/(I^3:x_1x_2x_3)\geq \left\lceil{\frac{\ell -5}{3}}\right\rceil.$  \QED

\begin{lemma}\label{CubeColonTwo}
Let $G$ be a graph and let $I=I(G)$. Fix $u \in V(G)$ and suppose that $xy \in E(G)$ with $x \in X^{\ell}$, where $X^{\ell}=X_G^{\ell}(u)$ for some $0 \leq \ell \leq d$. Then $\depth R/(I^3:xy) \geq \left\lceil{\frac{\ell -6}{3}}\right\rceil$.
\end{lemma}

\proof
We may assume $\ell \geq 7$ since otherwise the bound is trivial. First suppose either $x$ or $y$ is a leaf of $G$. Then by Lemma~\ref{leaf} we have that  $(I^3:xy)=I^2$ and by  Theorem~\ref{graphSquare}, we obtain $\depth R/(I^3:xy)\geq \left\lceil{\frac{d(I)-3}{3}}\right\rceil+p(I)-1$. Since $d(I) \geq \ell$, the result follows.

Next we assume that neither $x$ nor $y$ is a leaf of $G$. Let $\{z_1, \ldots ,z_s\}=N(x)\setminus\{y\}$ be ordered as in Lemma~\ref{OrderNeighbors}. Then, since $x,y,z_i \in \bigcup \limits_{j=\ell -1}^d X^j$, $\depth R_{i-1}/(I_{i-1}^3:xyz_i)\geq \left\lceil{\frac{\ell-6}{3}}\right\rceil$ by Lemma~\ref{CubeColonThree}. Now $x$ is a leaf of $I_s$, so $I_s^3:xy=I_s^2$, by Lemma~\ref{leaf}. Let $d(I_s)=d({}_{u}G(I_s))$. Then since $z_i \in \bigcup \limits_{j=\ell -1}^dX^{j}$, we have $d(I_s) \geq \ell-2$. Thus $\depth R_{s}/(I_s^3:xy) = \depth R_s/I_s^2 \geq \left\lceil{\frac{\ell-5}{3}}\right\rceil$ by Theorem~\ref{graphSquare}. Hence by Lemma~\ref{ExhaustTheNeighbors} we have $\depth R/(I^3:xy) \geq \left\lceil{\frac{\ell-6}{3}}\right\rceil$.
\QED

We are now ready to establish a bound for the depth of the third power of any edge ideal. 

\begin{theorem} \label{graphCube}
Let $G$ be a graph with $p$ connected components, $I=I(G)$, and let $d=d(G)$ be the diameter of $G$. Then $\depth R/I^3 \geq  \left\lceil {\frac{d-7}{3}} \right\rceil +p-1$. \end{theorem}

\proof
We proceed by induction on $n$, the number of vertices. We first handle the case when $n\leq 8$, in which case $d\leq 7$ and $p\leq 4$. When $p=4$, the result follows from Proposition~\ref{depth conn. comp} or from \cite[Theorem~3.4]{Mor}. By \cite[Lemma 2.1]{AJ} we know $\depth(R/I^t) \geq 1$ for all $t\leq p$, so the result holds for $p=3$. If $p=1$, or $p=2$ and $d \leq 4$, the bound is trivial. If $p=2$ and $d=5$, then the graph must consist of two disconnected paths, so the result follows from \cite[Theorem~3.4]{Mor}. Thus we may assume $n \geq 9$.

Let $u,v$ be the endpoints of a path that realizes the diameter of $G$ and let $X^i$ be as in Notation~\ref{XiSets}. Let $w \in N(v) \cap X^{d-1}$. 

Notice that $(I^3,w)=(J^3,w)$, where $J$ is the minor of $I$ formed by deleting $w$. We have two cases to consider. If $u$ and $v$ are in the same connected component of $J$ then $d(J) \geq d$ and $p(J)\geq p$, where $p(J)$ is the number of connected components of the graph associated to $J$. Hence by induction on $n$ we have $$\depth R/(I^3,w) \geq \left\lceil {\frac{d(J)-7}{3}} \right\rceil+p(J)-1 \geq \left\lceil {\frac{d-7}{3}} \right\rceil+p-1.$$ If $u$ and $v$ are not connected in $J$, 
 then $d(J) \geq d({}_{u}G(J)) \geq d-2$ and $p(J) \geq p+1$, or if $v$ is isolated, Lemma~\ref{addVariable} applies. Hence again by induction on $n$ we have 
$$\depth R/(I^3,w) \geq \left\lceil{\frac{d(J)-7}{3}}\right\rceil+p+1-1\geq \left\lceil{\frac{d-9}{3}}\right\rceil+p+1-1 \geq \left\lceil{\frac{d-7}{3}}\right\rceil+p-1.$$

Let $\{z_1, \ldots ,z_s\}=N(w)$ be ordered as in Lemma~\ref{OrderNeighbors}. Since $w \in X^{d-1}$ then by Lemma~\ref{CubeColonTwo} $\depth R_{i-1}/(I_{i-1}^3:wz_i)\geq \left\lceil{\frac{d-7}{3}}\right\rceil$.

Now $w$ is isolated in $I_{s}$ and thus $(I_s^3:w)=I_s^3$. 
Therefore by induction on $n$ we have that 
\begin{eqnarray*}
\depth R_s/(I_s^3:w)&=&\depth R_s/I_s^3 \geq \left\lceil{\frac{d(I_s)-7}{3}}\right\rceil +p(I_s)-1+1\\
&\geq& \left\lceil{\frac{d-3-7}{3}}\right\rceil+p-1+1 = \left\lceil{\frac{d-7}{3}}\right\rceil +p-1,\\
 \end{eqnarray*} since $d(I_s) \geq d({}_{u}G(I_s)) \geq d-3$ and $w$ is an isolated vertex.
Hence by Lemma~\ref{ExhaustTheNeighbors} we have that $\depth R/(I^3:w) \geq \left\lceil{\frac{d-7}{3}}\right\rceil+p-1$.

By applying the Depth Lemma~\cite[Proposition 1.2.9]{BH} to the following exact sequence $$0 \rar R/(I^3:w) \rar R/I^3 \rar R/(I^3, w) \rar 0$$
we have that $\depth R/I^3 \geq  \left\lceil {\frac{d-7}{3}} \right\rceil +p-1$.  \QED

As in Remark~\ref{remarkdistance} one may take $u$ and $v$ in the proof of Theorem~\ref{graphCube} to be endpoints of a path of length $\ell=d(u,v)$ and obtain $\depth R/I^3 \geq  \left\lceil {\frac{\ell-7}{3}} \right\rceil +p-1$. The next corollary follows from the proof of Theorem~\ref{graphCube}.

\begin{corollary}\label{CubeColonWithOne}
Let $G$ be a graph 
and let $I=I(G)$.  Fix $u \in V(G)$ and let $w \in X^{\ell}$ for some $0 \leq \ell $, where $X^i=X^i_G(u)$. Then $\depth R/(I^3:w) \geq \left\lceil{\frac{\ell -6}{3}}\right\rceil $. 
\end{corollary}

\proof  We may assume $\ell \geq 7$ since otherwise the bound is trivial.
Let $\{z_1, \ldots ,z_s\}=N(w)$ be ordered as in Lemma~\ref{OrderNeighbors}. By Lemma~\ref{CubeColonTwo} we have $\depth R_{i-1}/(I_{i-1}^3:wz_i)\geq \left\lceil{\frac{\ell -6}{3}}\right\rceil$.

Now $w$ is isolated in $I_{s}$ and thus $(I_s^3:w)=I_s^3$ and $d(I_s) \geq \ell -2$. Therefore by Theorem~\ref{graphCube}, we obtain 
$$\depth R_s/(I_s^3:w)\geq  \left\lceil{\frac{d(I_s) -7}{3}}\right\rceil +1 \geq \left\lceil{\frac{\ell -9}{3}}\right\rceil +1 =  \left\lceil{\frac{\ell-6}{3}}\right\rceil.$$ Hence by Lemma~\ref{ExhaustTheNeighbors}, the result follows.
\QED

The next example shows that the bound for the depth of the third power of an edge ideal given in Theorem~\ref{graphCube} is attained. This example extends naturally, which suggests a lower bound for the depth of any power.

\begin{example}\label{cube sharp ex}
{\rm
Let $R=k[x_1, \ldots, x_{10}]$ and let $I$ be the edge ideal of the graph $G$ below

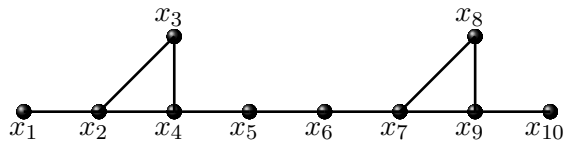
\begin{center}
\begin{tikzpicture}
\shade [shading=ball, ball color=black] (-1,0) circle (0.1) node [below ] {\textsf{\textbf{\textcolor{black}{\small{$x_1$}}}}};
\shade [shading=ball, ball color=black] (0,0) circle (0.1) node [below ] {\textsf{\textbf{\textcolor{black}{\small{$x_2$ }}}}};
\shade [shading=ball, ball color=black] (1,1) circle (0.1) node [above ] {\textsf{\textbf{\textcolor{black}{\small{$x_3$ }}}}};
\shade [shading=ball, ball color=black] (1,0) circle (0.1) node [below ] {\textsf{\textbf{\textcolor{black}{\small{$x_4$ }}}}};
\shade [shading=ball, ball color=black] (2,0) circle (0.1) node [below ] {\textsf{\textbf{\textcolor{black}{\small{$x_5$ }}}}};
\shade [shading=ball, ball color=black] (3,0) circle (0.1) node [below ] {\textsf{\textbf{\textcolor{black}{\small{$x_6$ }}}}};
\shade [shading=ball, ball color=black] (4,0) circle (0.1) node [below ] {\textsf{\textbf{\textcolor{black}{\small{$x_7$ }}}}};
\shade [shading=ball, ball color=black] (5,1) circle (0.1) node [above ] {\textsf{\textbf{\textcolor{black}{\small{$x_8$ }}}}};
\shade [shading=ball, ball color=black] (5,0) circle (0.1) node [below ] {\textsf{\textbf{\textcolor{black}{\small{$x_9$ }}}}};
\shade [shading=ball, ball color=black] (6,0) circle (0.1) node [below ] {\textsf{\textbf{\textcolor{black}{\small{$x_{10}$ }}}}};

\draw [line width=1pt  ] (-1,0)--(0,0);
\draw [line width=1pt  ] (0,0)--(1,1);
\draw [line width=1pt  ] (0,0)--(1,0);
\draw [line width=1pt  ] (1,0)--(1,1);
\draw [line width=1pt  ] (1,0)--(2,0);
\draw [line width=1pt  ] (2,0)--(3,0);
\draw [line width=1pt  ] (3,0)--(4,0);
\draw [line width=1pt  ] (4,0)--(5,0);
\draw [line width=1pt  ] (4,0)--(5,1);
\draw [line width=1pt  ] (5,0)--(5,1);
\draw [line width=1pt  ] (5,0)--(6,0);

\end{tikzpicture}

\end{center}

Then $d(G)=7$ and using Macaulay~2 \cite{M2} we have that $\depth R/I= \lceil{\frac{d+1}{3}}\rceil=2$, $\depth R/I^2= \lceil{\frac{d-3}{3}}\rceil=1$, and $\depth R/I^3= \lceil{\frac{d-7}{3}}\rceil=0$. Therefore, the bound in Theorem~\ref{graphCube} is sharp. 

Notice this is a graph with two copies of the graph in Example~\ref{square sharp ex} attached by an additional edge. One may attach more copies of the graph in Example~\ref{square sharp ex} to obtain examples of graphs where $\depth {R/I^t}=\left\lceil {\frac{d-4t+5}{3}} \right\rceil +p-1$ for any $t \geq 1$.  

}
\end{example}

Example~\ref{cube sharp ex} and Theorem~\ref{graphCube} lead to the following natural question.

\begin{question}\label{generalization}
Let $G$ be a graph with $p$ connected components, $I=I(G)$, and let $d=d(G) \geq 1$ be the diameter of $G$. Then is it true that for all $t\geq 1$ we have that $\depth {R/I^t} \geq \left\lceil {\frac{d-4t+5}{3}} \right\rceil +p-1$ or equivalently ${\rm{projdim}}_{R} R/I^t \leq n-\left\lceil{\frac{d-4t+5}{3}}\right\rceil -p+1$?

\end{question}

Clearly Theorems~\ref{graph}, \ref{graphSquare} and~\ref{graphCube} show that Question~\ref{generalization} has a positive answer for $t \leq 3$. If $t=4$ and $d\leq 2$, the bound in Question~\ref{generalization} reduces to $\depth R/I^t \geq p-4$ and so the result holds by Proposition~\ref{depth conn. comp}. Indeed, if $d \leq 2$, Proposition~\ref{depth conn. comp} gives a positive answer to the question for all values of $t$ and $p$. If $d=3$ and $t=4$, the bound reduces to $p-3$ and so is trivially true for $p \leq 3$. However, for any ideal $J$, $\depth(R/J) \geq 1$ if and only if $\m \not\in\Ass(R/J)$. By \cite[Lemma 2.1]{AJ} we know $\depth(R/I^t) \geq 1$ for all $t\leq p$, and so the question again has an affirmative answer for $p=4$. Thus the first case for which an answer to Question~\ref{generalization} is not known is when $d=3, t=4$, and  $p=5$.

An answer to this question would provide higher power analogues for Theorems~\ref{graphSquare} and~\ref{graphCube}. The difficulty in extending the outline of the proofs of those theorems to higher powers lies in generalizing the technical lemmas. For small powers, bounding the depth of $R/(I^t:x_1\cdots x_i)$ when the induced graph on $x_1, \ldots, x_i$ is connected is manageable because the small number of $x_i$ needed restricts the possible forms the induced graph can take. However, for higher powers of $t$, the products produced by repeatedly exhausting neighbor sets can induce graphs with poor behavior, including graphs for which $x_1\cdots x_i \not\in I^s$ for $s= \left\lfloor {\frac{i}{2}} \right\rfloor$.

We conclude this article by considering a few brief applications of our results. 
When $I$ is a square free monomial ideal then ${\rm projdim}R/I={\rm reg}(I^{\vee})$, where $I^{\vee}$ is the Alexander dual of $I$, \cite[Corollary~0.3]{Ter}. Since $I^{\vee \vee}=I$ then ${\rm reg}(I)={\rm projdim}(R/I^{\vee})=n-\depth R/I^{\vee}$, where $n=\dim R$. Using our result for the depth of the first power of edge ideals we may obtain bounds on these invariants as well. Another interesting invariant is Stanley depth. 
As a final application of our results we obtain lower bounds on the Stanley depth of the first three powers of edge ideals.

Let $R=k[x_1, \ldots, x_n]$ be a polynomial ring over a field $k$. 
Let $M$ be a nonzero finitely generated $\mathbb{Z}^n$-graded $R$-module, let
$u \in M$ be a homogeneous element and let $Z \subset \{x_1, \ldots, x_n\}$. Then $uk[Z]$ is the $k$-subspace generated by all monomials $uv$, where $v$ is a monomial in $k[Z]$. A presentation of $M$ as a finite direct sum of such spaces $\mathcal{D}$: $M=\bigoplus \limits_{i=1}^{r}u_ik[Z_i]$ is called a {\it{Stanley decomposition}} of $M$. Let $\sdepth {\mathcal{D}}=\min\{|Z_i|: i=1, \ldots, r\}$ and let $\sdepth M=\max\{ \sdepth{\mathcal{D}}: \mathcal{D} \mbox{ is a Stanley decomposition of } M\}$. Then $\sdepth{M}$ is called {\it{Stanley depth}} of $M$.  It was conjectured by Stanley in \cite{Sta} that $\sdepth{M} \geq \depth M$ for all $\mathbb{Z}^n$-graded modules $M$. There has been considerable interest  concerning this conjecture, see for instance \cite{Her}, and recently  Duval, Goeckner, Klivans, and Martin found a counterexample to Stanley's conjecture,  \cite{DGKM}. However, finding classes for which the conjecture holds is still an interesting endeavor. 

For the case of edge ideals of graphs and their powers we are able to obtain lower bounds for the Stanley depth using our results from the previous sections as well as the following version of the Depth Lemma for Stanley depth.

\begin{lemma}\cite[Proposition~2.6]{BKU}, \cite[Lemma~2.2]{Ra} \label{sdepth lemma}
Let $R=k[x_1, \ldots, x_n]$ be a polynomial ring over a field $k$. Let $0 \rar M \rar N \rar L \rar 0$ be a short exact sequence of finitely generated $\mathbb{Z}^n$-graded $R$-modules. Then $\sdepth{N} \geq \min\{\sdepth{M}, \sdepth{N}\}$.

\end{lemma}

\begin{theorem}\label{Stanley}
Let $G$ be a graph with $p$ connected components,  $I=I(G)$, and let $d=d(G)$ be the diameter of $G$. Then for $1 \leq t \leq 3$ we have
$$
\sdepth {R/I^t} \geq \left\lceil {\frac{d-4t+5}{3}} \right\rceil +p-1.
$$
\end{theorem}

\proof
The proof follows by induction on $n$, the number of vertices of $G$. Given Lemma~\ref{sdepth lemma} we can proceed the same way as in the proofs of Theorems~\ref{graph},~\ref{graphSquare},~\ref{graphCube} as long as we can establish the bounds for the base case of the induction, that is when $n=d+1$ and $G$ is the graph of a path. The required bounds are known to hold for the Stanley depth, see for example \cite[Theorem~2.7]{PFY}. \QED

One consequence of Theorem~\ref{Stanley} is that any class of ideals for which at least one of the bounds in Theorems ~\ref{graph},~\ref{graphSquare},~\ref{graphCube} is an equality will correspond to a class of modules that satisfy the Stanley conjecture. Thus discovering when the bounds are achieved is an area of further interest.

\section{Acknowledgments}
The authors would like to thank  Nate Dean for helpful conversations regarding graph theory and the anonymous referees for useful comments and suggestions. The first author would also like to thank the Department of Mathematics at Texas State University for its hospitality while some of the work was completed.

\end{document}